\documentclass[10pt,twocolumn,twoside]{IEEEtran}
\newcommand{\comm}[1]{}
\usepackage{amssymb}
\usepackage{amsmath}
\usepackage{mathrsfs}

\usepackage{epsfig}

\def\xxxonly{\comm}
\def\noxxx{\comm}
\def\blind{\comm}\def\blind{ }
\def\xxxonly{ }
\usepackage{xcolor}

\definecolor{my}{rgb}{0.1, 0.5, 0.4}

\newtheorem{theorem}{Theorem}[section]

\def\e{\varepsilon}

\def\defi{\stackrel{{\scriptscriptstyle \RR}}{=}}

\def\a{\alpha}
\def\d{\delta}

\def\F{{\cal F}}
\def\w{\widehat}
\def\Ind{{\,\rm Ind\,}}
\def\Ind{{\mathbb{I}}}

\def\const{{\rm const\,}}

\def\R{{\bf R}}
\def\E{{\bf E}}
\def\P{{\bf P}}

\def\L{L}

\def\b{\beta}
\def\s{\delta}

\def\W{{\cal W}}
\def\ww{\widetilde}

\def\t{\theta}
\def\oo{\bar}
\def\s{\sigma}

\def\p{\partial}
\def\G{\Gamma}

\def\V{{\cal V}}

\def\L{{\cal L}}

\def\LL{{\oo{\bf L}}}

\newcommand{\be}{\begin{equation}}
\newcommand{\ee}{\end{equation}}
\newcommand{\bd}{\begin{displaymath}}
\newcommand{\ed}{\end{displaymath}}
\newcommand{\ba}{\begin{array}{ll}}
\newcommand{\ea}{\end{array}}
\newcommand{\baa}{\begin{eqnarray}}
\newcommand{\eaa}{\end{eqnarray}}
\newcommand{\baaa}{\begin{eqnarray*}}
\newcommand{\eaaa}{\end{eqnarray*}}
\font\sm=cmr10

\def\oo{\bar}

\usepackage{graphicx}

\def\LL{\mathscr{L}}
\def\UU{{\mathscr{U}}}
\def\RR{R}
\title{
\noxxx{Optimal energy storing and dispatching with  preferable multi-battery charging regimes}
\xxxonly{Optimal energy storing and selling  in continuous time stochastic multi-battery setting}} 
\author{
Nikolai Dokuchaev}
\begin{document}
\def\brea{}
\def\breakk{}
\def\break{}
\def\break{\nonumber\\ }\def\breakk{\nonumber\\&&}\def\brea{\nonumber\\}\def\BR{\vspace{-1mm}\\ & }

\maketitle
\comm{\let\thefootnote\relax\footnote{  ITMO University, St. Petersburg, 197101 Russian Federation }}

\begin{abstract}
The paper  suggests a  new stochastic model for energy producing, dispatching, and storing
 in the multi-battery setting
that takes into account the topology of the system of the links between the batteries, the transmission and storage losses, and  requirements for special regimes for batteries charging and discharging
helping to prolong batteries life. For this model, the problem of optimal energy storing and dispatching
is considered and solved using dynamic programming and duality methods.

\xxxonly{\par{\bf 2010 Mathematics Subject Classification}: 91B70,  
93E20
}\end{abstract}

\section{Introduction}
Recent widespread expansion of new energy technologies and growth in the number of
small and medium  companies producing and selling energy has triggered need for new types of operating algorithms
 to ensure sustainable and reliable production process and  maximization of the profit.
 An important feature of energy production based on the renewable sources is that unpredictable fluctuations of the production rate can be significant.
 To compensate these fluctuation and to ensure a
 more stable output level,  the energy must be stored. Typically, it is necessary to consider a storage consisting   of several separate battery  units that have to be regularly charged and discharged.

 The fluctuations of the production rate  represent a mixture of relatively regular predictable components such as night interruptions for solar energy
 and tide cycles for wave energy, and of irregular unpredictable components such as fluctuations caused by weather conditions for solar and wind energy;
 see e.g. \cite{U1,U2,U3}.    The fluctuations  of the selling price
 also represent a mixture of relatively regular  and predictable components such  as switching between day and night prices and irregular unpredictable components caused  by
 unpredictable market movements.   The production  rate and selling prices
 may have unpredictable deviations; their   forecast
 without error is impossible.

Therefore, the  power producers need optimal strategies for energy storing and selling  that help to reduce the impact of unpredictability of production rate and market prices.  The problem appears to be a  control problem under uncertainty, that, in the case of renewable energy,  is defined by non-controlled external factors such as the
weather.  An important feature of  this problem is that the dimension of the control process may be high for systems with many production and storage units. 
These questions are important for applications in energy sector, particularly for small and medium producers of renewable energy.
 The problem
was studied intensively; see e.g. \cite{1,CM, Kim,O,A0,A1,A2,W}.

The present paper suggests a comprehensive and yet compact dynamic model of energy dispatching  and storage  for a  multi-battery setting.
This model represents a  development of  models suggested in \cite{1,A0,A1,A2,W}.
The novelty is that the topology of the links between the batteries and transmission/storage  are taken into account. Another  novelty is that our setting takes into account  a special feature of energy trading:
the storage is based on batteries requiring  certain regimes of charging and discharging to prolong  the battery life; see, e.g., \cite{Deg,Ning,A0} and the bibliography therein.
Given that the batteries are expensive, this is  a significant factor in decision making.
To address this, we
considered an extended model where cumulative moving averages were included.

 The main focus of the paper is modelling of the control problem for energy storing and dispatching. In addition, the paper suggests some approaches  for the solution of the corresponding
 stochastic optimal control problem.
 For this control problem, the objective is to select the regimes for supplying (selling) the energy to the grid, the
depositing the energy in the batteries, and for redistributing the energy among the batteries, with a performance
criterion that takes into account the obtained monetary gain and the regimes for the batteries.

The main feature of
this problem is that  there is a fixed  domain where
 state processes are  allowed to go on and off the  boundary of the admissible domain as well as stay on the boundary. In known  stochastic optimal control solutions,
problems with boundary are usually  considered  for
processes with reflection from the boundary or with absorption on the boundary.

For the setting with
Markov diffusion model for the random factors, the paper suggests a solution
based on the dynamic programming method. \comm{following  the  approach used for  energy trading models in \cite{D13}.}
We derived the equation for the optimal value of the problem in a form of a Hamilton-Jacobi-Bellman (HJB) equation,
and obtained some existence results.
  For a large number of factors arising in a multi-battery setting, the state space dimension for the HJB equation could be high, and the  numerical solution could be challenging. To address this, the paper suggests  an  approach based on duality and pathwise optimization, in the spirit of  \cite{BS,BD2,Br,Davis,HK}.
In the framework of this approach, the optimal value function can be calculated using Monte-Carlo simulation of the Lagrange terms and pathwise deterministic
optimization in a class of non-adapted processes. This does not lead to an optimal strategy immediately; however, it
gives an opportunity to estimate how far from optimal is the performance of a particular  strategy, for instance, such as suggested in
\cite{1}.

 The paper is organised as follows: Section \ref{sec1} describes the basic model setting with a single battery.   Section \ref{secM} introduces a multi-battery setting and discuss optimization of battery regimes. Section \ref{secHJB} discusses   Hamilton-Jacobi-Bellman equations
 describing the optimal value functions.  Section \ref{secD} suggests a duality and pathwise optimization approach   for estimation of the optimal value functions.
 Section \ref{secP} contains the proofs. Section \ref{secC} offers some discussion and concluding remarks.

\section{Problem setting for the basic model}\label{sec1}
In this section, we present a  simplified  version of the model to outline some basic features.
We consider a model consisting of an energy producing plant  and a battery storage operating under the common management.
 This model represents a Virtual Power Plant: its purpose is to  supply the energy in  the external grid.

 Assume that $p(t)$ is a random
process representing the current rate of production of the energy by the plant, and that $S(t)$ is a random
process representing the current price of  the energy unit at time  $t\in[0,T]$, where $T>0$ is a given terminal time. These processes may depend on unpredictable
factors (weather, market conditions).

The objective of the controller  is to select the regimes of supplying (selling) the energy to the grid, and of
depositing the energy in the battery. More precisely, we assume  that the controller
 has to   select   the rate  $u(t)$ of depositing  the energy   into the storage (battery).
This also defines the rate $q(t)$ of the supplying energy to the external grid as
 \baaa
 q(t)=p(t)-u(t).\eaaa
The process $u(t)$ can be considered a control (strategy).

 The case where $q(t)<0$ is not excluded; this would correspond to the case where  the plant withdraws (buys) the energy from the grid and  deposit it in its battery.

Let $y(t)$ represents the quantity of the energy currently stored in the battery
such that  \baaa
\frac{dy}{dt} (t)=u(t)-\a y(t),
\label{yC}\eaaa
 where $\a>0$ is a parameter representing the storage losses.

We assume that    \baaa
u(t)\in [-L,L],\qquad y(t)\in[0,C],
\label{uy}\eaaa
where $L>0$ is a parameter describing  restrictions on the energy transfer rate, and $C>0$
is a parameter representing  the battery capacity.

The control process $u(t)$  has to be selected using the historical
observations of $(p(t),S(t))$  as  well as other currently available information such as  weather data or currency exchange rate.

Let $T>0$ be a given terminal time. We assume that the monetary value of the output of a particular strategy $u$ for this $T$
can be represented as \baaa
&&F(u)=\int_0^Tq(t)S(t)dt+  S(T)y(T)\breakk=\int_0^T(p(t)-u(t))S(t)dt+  S(T)y(T) .\eaaa
The integral part here is the  value representing the total earning from the selling during the time period $[0,T]$. The value  $S(T)y(T)$ represents the
market value  of the stored energy at the terminal time $T$. The goal is to maximize  $F(u)$ over $u$.

The paper focuses on the setting where  the future values of $(p(t),S(t))$ are random, and their future values to be forecasting
with possible forecasting error. In this case, the goal is to maximize the expectation $\E F(u)$ over $u$ given a probability distribution describing the current
hypothesis on $(p,S)$ and on other factors.

Let us give a more accurate  description of information available for the decision making.

Let $\{\F_t\}_{t\in[0,T]}$
be the filtration representing the information given the current and past observations available at time $t$. The processes  $p$, $S$, $q$, and $u$ have to be $\F_t$-adapted.

 This filtration may also include information generated by other processes such as the weather etc.

\par

We consider strategies   $u(t)$ that are  $\F_t$-adapted  and such that $-L\le u(t)\le L$ for all $t$.

 The following stochastic optimal control problem arises: for given $t<T$ and $\eta\in [0,C]$,
 \baa
&&\hbox{Maximize}\quad\breakk\E\left[\int_t^T(p(s)-u(s))S(s)ds +\ y(T)S(T)\right] \quad\breakk\hbox{over}\quad
u\nonumber
\\ &&\hbox{subject to}\quad  \frac{d y}{ds}(s)=u(s)-\a y(s),\quad\nonumber
\\&& y(t)=\eta,  \quad y(s)\in[0,C]\quad\hbox{for}\quad s\in[t,T].
 \label{O1}\eaa
Here $\E$ is the expectation.
 \par
It can be noted that admissible
 state processes  $y$  are  allowed to go on and off the  boundary of the admissible domain  as well as stay on the boundary. This setting
is non-standard  for stochastic optimal control theory,
where it is more common to consider processes with reflection from the boundary or with absorption on the boundary.
\begin{theorem}\label{ThD1}
Problem (\ref{O1}) is equivalent to the problem
 \baa
&&\hbox{Maximize}\quad\breakk\E\left[\int_t^T(p(s)-u(s))S(s)\Ind_{\{0\le y(s)\ge C\}}ds + y(T)S(T)\right] \quad\breakk
\breakk\hbox{over}\quad
u\nonumber \\&&
\hbox{subject to}\quad  \frac{d y}{ds}(s)=(u(s)-\a y(t))\Ind_{\{0\le y(s)\le C\}},\quad\breakk y(t)=\eta.
 \label{K1}\eaa
\end{theorem}
Here $\Ind$ is the indicator function.

The cases where  $p(t)<0$, $q(t)<0$, or $u(t)<0$, are not excluded: in this case,  $p(t)$ is the rate of losing energy (this could occur, for instance, due to
  technological issues),  $q(t)<0$ is the rate of buying the energy, and  $u(t)<0$ is the rate of withdrawing the energy  from the storage.
 The case of a non-positive  $S(t)$  is also  not excluded, even if this is  a rare possibility (in calculations, this possibility can be taken into account via an appropriate choice  of a model for the energy prices).
\subsection{The case of stochastic Markov model}
Up to the end of Section \ref{sec1}, we consider a stochastic model  for the process $(p(t),S(t))$. We assume that it evolves as  a part of stochastic Markov
diffusion process (see, e.g., \cite{Krylov}).

Let  $w(\cdot)$ be a standard $n$-dimensional Brownian motion process, $n\ge 2$ such that
$\E w(t)=0$ and $\E w(t)^2=t$, $t>0$. Let
 $\oo g:\R^n\times[0,T]\to\R$ and
$\oo\b:\R^n\times[0,T]\to\R$ be some continuous functions such that
$|\oo g(x,t)|+|\oo \s(x,t)|\le \const(|x|+1)$ and  $|\p \oo g(x,t)/\p x|+|\p \oo \s(x,u,t)/\p x|\le \const$ for all $x,u,t$.

In this section, we assume that $\{\F_t\}$ is the filtration generated by $w(t)$,
\baa
 &&p(t)= x_1(t),\quad S(t)=x_2(t), \label{diff}
\eaa
where $\oo x(t)=(\oo x_1(t),...,\oo x_n(t))^\top$  is a stochastic
diffusion process evolving as
\baa
&& d\oo x(t)=\oo g(\oo x(t),t)dt+ \oo\s(\oo x(t),t)dw(t).
\label{oox}
\eaa
The components$\{\oo x_k(t)\}_{k>2}$  represent currently available but unpredictable information  (other  than $(p(t),S(t))$)
 such as the weather data or a currency exchange rate.

Equations (\ref{diff})-(\ref{oox}) define a stochastic evolution model for the process $(p(t),S(t))$.
Calibration of  the parameters for these
 equations  is a complicated task involving statistical inference and  forecasting methods; this problem
 was considered e.g. in \cite{S2}.

Matching of the definitions shows that  problem (\ref{K1})  can be rewritten in the form of  a problem
 \baa
&&\hbox{Maximize}\quad\breakk\E\left[\int_s^T h(x(s),u(s),s)ds  +\Phi(x(T))\right] \quad\hbox{over}\quad
u\nonumber\\ &&
\hbox{subject to}\quad \breakk  d x(s)=g( x(s),u(s),s)ds+\s(x(s),s)dw(s), \quad\breakk x(t)=\xi.
 \label{K}\eaa
\renewcommand{\arraystretch}{1.2}
Here \baaa
&&\xi=\left(\begin{array}{c}
                        \oo x(t) \\
                        y(t)
                      \end{array}\right)\in\R^{n+1},\quad x(s)=\left(\begin{array}{c}
                        \oo x(s) \\
                        y(s)
                      \end{array}\right),\quad
\breakk g(x,u,s)=\left(\begin{array}{c} \oo g(x_1,...,x_n,s)\\  f(x_{n+1},u,s) \end{array}\right),
\quad\\ &&\s(x,u,s)=\left(\begin{array}{c}\oo \s(x_1,...,x_n,s)\\0_{\R^{1\times n}}\end{array}\right), \quad
\breakk f(x_{n+1},u,s)= (u-\a x_{n+1})\Ind_{\{0\le x_{n+1}\le C\}},\quad\\&&
h(x,u,s)=(x_1-u)x_2\Ind_{\{0\le x_{n+1}\le C \}},\quad
\breakk \Phi(x)=x_2 x_{n+1},
\qquad
\eaaa where $x=(x_1,...,x_{n+1})\in\R^{n+1}$, $s\in[t,T]$.
 \section{Multi-battery model with optimization of the battery regimes}\label{secM}
 Consider now situation where the storage consists  of $m$ separate but linked batteries
 with different regimes of their operations, where $m\ge 1$.

We assume that the controller selects the rate of selling  energy to the external grid and the rate of depositing  energy from the plant  into each particular battery.
In addition, the controller  selects the rate of energy transfers between the batteries.

 In other words, the controller has to calculate  vector processes  $u(t)=\{u_i(t)\}_{i=1}^m$  and matrix processes
 $v(t)=\{v_{ij}(t)\}_{i,j=1}^m$, where
  $u_i(t)$ represents the rate of depositing   the energy   into the battery $i$ from the plant, and $v_{ij}(t)$
is the rate of energy transferred from the battery $i$ to be deposed in the battery $j$.

Let $y_i(t)$ be the quantity of the energy stored in the  $i$th battery, and let $y(t)=\{y_i(t)\}_{i=1}^m$.

 To take this into account, we
extend the model introduced above.

Similarly to the case of a single battery considered above, the rate $q(t)$ of selling the energy to the external grid can be represented as
 \baaa
 q(t)=p(t)-\sum_{i=1}^m u_i(t).\eaaa
Here $p(t)$ is ahain the rate of production  of energy.

For the dynamics of the storage levels, we develop below a more advanced model
that takes
into account the topology of the links, storage  losses, and link losses.

\subsection{The class of admissible strategies}

Let $L_{i}\ge 0$ be given for $i=1,...,m$, and let
$\oo L_{ij}\ge 0$ and $\ww L_{ij}\ge 0$ be given for $i,j=1,...,m$
such that \baaa
\oo L_{ij}=-\ww L_{ji},\quad \ww L_{ij}=-\oo L_{ji}.
\eaaa

Let sets $\RR({\rm p,y})\subset \R^m\times \R^{m\times m}$  be defined for ${\rm p}\in\R$ and ${\rm y}\in\R^m$.
\par
Let $\UU$ be the class of pairs $(u,v)$  such that $u(t)$ and $v(t)$ are $\F_t$-adapted and
\baa
&&u_{i}(t)\in [-L_{i},L_{i}],\quad v_{ij}(t)\in [-\oo L_{ij},\ww L_{ij}],\label{uv}
\\
&&  v_{ij}(t)\equiv  -v_{ji}(t).\label{vv}
\\
&&(u(t),v(t))\in \RR(p(t),y(t)),\quad t\in[0,T]\label{uvy}
\eaa

The choice of the sets $\RR({\rm p,y})$ allows to impose various  restrictions for the model. For example,
selection of  $\RR({\rm p},{\rm y})=\{{\rm u}\in\R^m, {\rm v}\in \R^{m\times m}:\ \sum_i {\rm u}_i\le p\}$
is appropriate for a model where the energy is not purchased from the external sources.
Another example: The choice of  $\RR({\rm p},{\rm y})=\{{\rm u}\in\R^m, {\rm v}\in \R^{m\times m}:\  {\rm v}_{ij}\le 0\ \hbox{if}\  {\rm y}_i\le
{\rm y}_j \}$
is appropriate for a model where the energy can be transferred from battery $i$ to battery $j$
only if $y_i(t)>y_j(t)$.
 \par

The case where $L_{i}=0$ is not excluded; this would mean that $u_{i}(t)\equiv 0$,
i.e., there is no energy transmission from the plant to the battery $i$.

The case where $\oo L_{ij}=\ww L_{ij}=0$ is also not excluded; this would mean that $v_{ij}(t)=v_{ji}(t)\equiv 0$,
i.e., there is no energy transmission between the batteries $i$ and $j$.

This implies  that the choice of the non-zero values
$L_i$, $\oo L_{ij}$, and $\ww L_{ij}$, defines the topology  of  the system plant/batteries, i.e., it defines the  links between  the batteries and the plant and
the mutual links between  the batteries.

For a model where the energy is not purchased from the external sources, one  could consider restriction that $u(t)\in [-L,\min(p(t),L)]$. We omit this case.

\subsection{The evolution of the batteries loads}

To take the links between batteries and links/storage losses into account, we accept the following  model:
 \baaa
&& \frac{dy_i}{dt} (t)=u_i(t)-\a_i y_i(t)-\sum_{j:\ j\neq i}\b_{ij}(v_{ij}(t)), \quad
\breakk y_i(t)\in[0,C_i].
 \eaaa
The parameters $\a_i\in [0,1)$ describe the storage losses for battery $i$; they may depend, for instance,
on the age of a battery. The parameters $C_i>0$ describe the  capacity of the battery $i$.

The functions $\b_{ij}:\R\to [0,1]$ are such that
$\b_{ij}(v)=v$ for $v\ge 0$ and $\b_{ij}(v)=b_{ij} v$ for $v< 0$, where $b_{ij}\in [0,1]$
are parameters describing the rate of the losses for transmission  from $i$ to $j$.
These parameters $b_{ij}$ may depend, for instance, on the distance between the batteries.
In particular, due the transmission losses, the
rate of energy  depositing  in the battery $j$ from the battery $i$ can be  less than than
the rate of withdrawing energy from battery $i$ for the battery $i$.
 For example,
if $b_{ji}=0.8$  and $v_{ij}(t)>0$ then $-\b_{ji}(v_{ji}(t))=-0.8 v_{ji}(t)=0.8 v_{ij}(t)>0$ is the rate of energy depositing in the battery $j$ from the battery $i$; on the other hand,  there is  energy withdrawal  from the battery $i$ for the battery $j$;
with the rate $-v_{ij}(t)=v_{ji}(t)<0$.

We consider below vector processes $y(t)=(y_1(t),...,y_m(t))$ and $u(t)=(u_1(t),...,u_m(t))$.
In addition, we consider matrix process $v(t)=\{v_{ij}(t)\}_{i,j=1}^m$.

 \subsection{Taking into account preferable regimes}

 The technology reasons suggest certain regimes for  charging and discharging batteries used to store energy by the producer.
 Given that the batteries are expensive, this could be a significant factor in decision making.
 It is known that charging and discharging too rapidly may lead to shortened battery life \cite{Deg}. This can be controlled by using $L$ in our setting.
 In addition, a deep discharge may also have  negative effect  \cite{Deg}. To take this into account, we may incorporate the additional objective of maximization of
   \baa
   \E \int_0^T \phi(y(t))dt,
   \label{wphi}\eaa
where $\phi:\R^m\times \R^m\to(-\infty,0)$ is a function
that achieves minimum on the boundary of the rectangle domain  $\prod_{i=1}^m [0,C_i]$ (or on a selected part of the boundary). For example, one
may select \baaa
\phi(y)=-\prod_{i=1}^m\varphi_i^y(y_i),\eaaa where  $\varphi_i^y:[0,C_i]\to(0,+\infty)$ are some
$U$-shaped convex functions.
\subsubsection*{Preferences using cumulative moving averages}
It appears that some important performance indicators  cannot be described by integrals of   functions of the current state $y(t)$.
In some cases, it could be reasonable to use performance criterions  involving
 integral functionals on the paths, such as cumulative moving averages \baaa
\oo y_i(t)\defi \frac{1}{t}\int_0^ty_i(s)ds.
\eaaa
This can be described as  via minimization
of the expectation
  \baa
   \E \int_0^T \phi(y(t),\oo y(t))dt,
   \label{ooy}\eaa
where $\oo\phi:\R^m\times \R^m\times(0,T]\to\R$ is a function describing the agent's preferences.

 For instance,  a preference that the charging processes are oscillating with a similar rate for all batteries can be taken into account with
  \baaa
   \phi(y(t),\oo y(t))= -\sum_{i,j=1}^m (Y_i(t)-Y_j(t))^2,
   \label{ooy1}\eaaa
  where $Y_i=y_i(t)-\oo y_i(t)$.
The corresponding performance indicator cannot be quantified via (\ref{wphi}).

Alternatively, one may prefer to have  all batteries have similar charges.
 For this, one can use
  \baaa
\phi(y(t),\oo y(t))= \sum_{i,j=1}^m(\oo y_i(t)-\oo y_j(t))^2.
  \eaaa
\par
Regularity of  cycles of particular batteries can be controlled via  maximization of (\ref{ooy}) with
  \baa
   \phi(y(t),\oo y(t))=-\G\sum_{i=1}^m(y_i(t)-\oo y_i(t))^2.
   \label{phiooy}\eaa
  Clearly, the maximum  of (\ref{ooy}) with $\G>0$  is achieved for the batteries with constant levels of energy stored.
The maximum of (\ref{ooy}) with $\G<0$  is achieved for the batteries with oscillating  energy levels.
\def\yy{{\rm{y}}}

\begin{figure}[ht]
\centerline{\includegraphics[width=9.5cm,height=6.5cm]{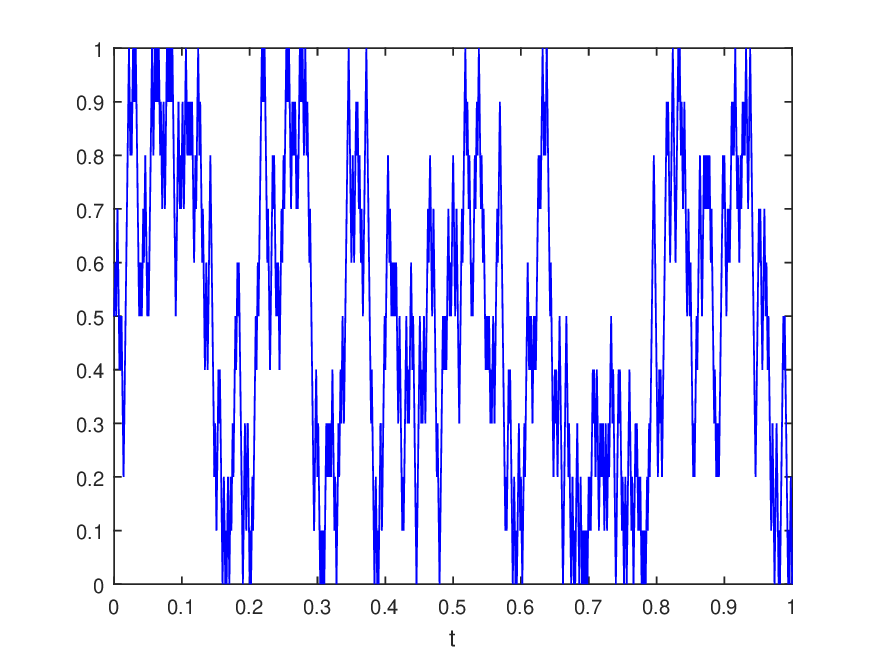}} 
\caption[]{\sm
Approximation of optimal solution of problem (\ref{bat})}
\label{fig}
\end{figure}

   Let us demonstrate the impact  of maximization of (\ref{ooy}) given
   (\ref{phiooy}) with $\G>0$. Let $t\in(0,T)$. Consider the following problem:
\baa
&&\hbox{Maximize}\quad \int_t^T(\yy(s)-\oo\yy(s))^2ds \quad \hbox{over}\quad \yy(\cdot),
\\
&&\hbox{subject to}\quad\breakk \yy(s)\in[0,C],\quad \frac{d\yy}{ds}(s)\in[-L,L],
\label{bat}
\eaa
where
\baaa
\oo\yy(s)\defi \frac{1}{s}\int_0^s\yy(\t)d\t.
\eaaa
The optimal $\yy(s)$  must deviate from their historical mean as much as possible.
In fact, problem (\ref{bat}) can represented as a linear quadratic problem
\baaa
&&
\hbox{Maximize}\quad \int_t^T(\yy(s)-\oo\yy(s))^2ds\quad \hbox{over}\quad \yy(\cdot),
\nonumber\\
&&\hbox{subject to}\quad \breakk d\oo\yy(s)/ds=-s^{-2}\oo\yy(s)+s^{-1}\yy(s),\quad\breakk
\yy(s)\in[0,C],\quad \frac{d\yy}{ds}(s)\in[-L,L].
\label{bat2}\eaaa
It appears that this criterion lead to periodic regimes with stable oscillations.
To show this, we did the following experiments.  We created a set of discrete time paths
$u=(u^{(1)},....,u^{(N)})$ such that $u^{(k)}=\pm C$, where $N$ is the time discretization parameter.
This set of paths was created using Monte-Carlo simulation of binary vectors with independent components.
The corresponding process $\yy(s)$  was replaced by the vector  $\yy=(\yy^{(1)},...,\yy^{(N)})$
such that    $\yy^{(j)}=\yy^{(1)}+ \sum_{d=1}^j\yy^{(d)}\RR t$, where $\RR t=T/N$. The approximation of the optimal path
was identified as the path  with the minimal value of
 \baaa
\sum_{j=1}^n (\yy^{(j)}-\oo\yy^{(j)})^2
\label{bat1}\eaaa
where
\baaa
\oo\yy^{(j)}= j^{-1}\sum_{d=1}^j \yy^{(d)}.
\eaaa
The approximation of the path among  Monte-Carlo simulated 2,000 paths is presented by Figure \ref{fig}.
For this experiment, we used $T=1$, $C=1$, $L=100$, $n=1000$.
\subsection{Optimal control setting for the multi-battery model}
\begin{table}[ht]
 \begin{tabular}{|ll|}\hline
 \multicolumn{2}{|c|}{{\bf The list of the model parameters and notations }}
 \\
 $p(t)$ & The production rate\\
 $S(t)$ & The current energy price
 \\
 $\{\F_t\}$ & The filtration representing the flow of \BR available information
 \\
 $m$   & The number of batteries \\
 $C_i$ & The capacity of the battery $i$\\
 $y_i(t)$ & The quantity of the energy currently stored \BR in the battery $i$\\
 $\oo y_i(t)$ & The cumulative moving average  of $y_i(t)$\\
 $u_i(t)$ & The controlled  rate of energy transfer to the  \BR battery $i$ from the production plant \\
   $[-L_i,L_i]$ & The interval of possible values of $u_i(t)$\\
 $v_{ij}(t)$ & The controlled rate of energy transfer from \BR the battery $i$ to the battery $j$ \\
 $[-\oo L_{ij},\ww L_{ij}]$ & The interval of possible values of $v_{ij}(t)$
 \\
 $\UU$ & The set of $\F_t$-adapted  $(u(t),v(t))$ such that \BR (\ref{uv})-(\ref{uvy}) holds. \\
 $\a_i$ & The rate of storage losses for the battery $i$\\
 $b_{ij}$ & The rate of transfer losses for the energy transfer \BR from the battery $i$\ to the battery $j$\\
$\b_{ij}(v)$ & $\b_{ij}(v)=v$ if $v>0$, and  $\b_{ij}(v)=b_{ij} v$ if $v>0$\\

 $\E$ & The expectation\\
 $ \Ind$ & The indicator function \\
\hline
 \end{tabular}
\label{tabA}\end{table}

Let a function $\phi:\R^n\times \R^n\times [0,T]$ be selected  with the purpose of taking into account the preferences for the battery regimes.

 The following stochastic optimal control problem arises: for given $t\in (0,T)$ and $\eta\in \prod_{i=1}^m[0,C_i]\subset\R^m$,
 \baa
&&\hbox{Maximize}\quad\nonumber\\
&&\E\int_t^T\left[\left(p(s)-\sum_{i=1}^m u_i(s)\right)S(s)+\phi(y(s),\oo y(s),s)\right] ds \breakk+\E S(T)\sum_{i=1}^m y_i(T)
\quad\nonumber\\&&\hbox{over}\quad
(u,v)\in \UU,\nonumber \\&&
\hbox{subject to}\quad  \frac{d y_i}{ds}(s)=u_i(s)-\a_i y(s)
-\sum_{j:\ j\neq i}\b_{ij}(v_{ij}(s)),\nonumber\\&&
\quad y_i(s)\in[0,C],\quad i=1,...,m,\quad y(t)=\eta.
 \label{BO1m}\eaa
 Here  $\oo y_i(s)\defi \frac{1}{s}\int_0^s y_i(q)dq$.
 \begin{theorem}\label{Th1M}
Problem (\ref{BO1m}) is equivalent to the problem
 \baa
&&\hbox{Maximize}\quad\nonumber\\
&&\E\int_t^T\Biggl[\left(p(s)-\sum_{i=1}^m u_i(s)\right)S(s)\Ind_{\{0\le y_{i}(s)\le C,\ i=1,...,m \}}\breakk +\phi(y(s),\oo y(s),s)\Biggr]ds+\E S(T)\sum_{i=1}^m y_i(T)\quad\nonumber
\\&&\hbox{over}\quad
(u,v)\in \UU,\nonumber \\&&
\hbox{subject to}\quad\nonumber \\&& \frac{d y_i}{ds}(s)=\Bigl(u_i(s)-\a_i y(s)\breakk -\sum_{j:\ j\neq i}\b_{ij}(v_{ij}(s))\Bigr)\Ind_{\{0\le y_i(s)\le C\}},\quad i=1,...,m,\nonumber\\&&\quad y(t)=\eta.\hphantom{5cm}
 \label{BKm}\eaa
\end{theorem}
It can be noted that we do not exclude cases where the process $(p(t),S(t))$.
In this case,  (\ref{BKm}) is  a deterministic optimal control problem.
\subsection{The case of stochastic Markov multi-battery model}\label{SubsecDif}
In this section, we sssume that $(p(t),S(t))$  is a stochastic Markov
diffusion process $\oo x(t)$  defined as (\ref{diff})-(\ref{oox}), and that
  $\{\F_t\}$ is generated by $w(t)$.
 Matching of the definitions shows that  problem (\ref{BKm}) is equivalent to the problem
 \baa
&&\hbox{Maximize}\quad\E\left[\int_t^T h(x(s),u(s),s)ds  +\Phi(x(T))\right] \quad\breakk\hbox{over}\quad
(u,v)\in \UU\nonumber \\
&&\hbox{subject to}\quad\breakk d x(s)=g(x(s),u(s),v(s),s)ds+\s(x(s),s)dw(s), \quad\breakk x(0)=\eta.
 \label{DKm}\eaa
\renewcommand{\arraystretch}{1.2}
Here $t\in[0,T]$,  $x(s)=(x_1(s),...,x_{n+2m}(s))^\top$,   $x_1(s)=\oo x_1(s)=p(s)$, $x_2(s)=\oo x_2(s)=S(s)$,
$(x_{n+1}(s),...,x_{n+m}(s))^\top=(y_{1}(s),...,y_{m}(s))^\top$, $(x_{n+m+1}(s),...,x_{n+2m}(s))^\top=(\oo y_{1}(s),...,\oo y_{m}(s))^\top$,
 $f=(f_1,...,f_{2m})^\top$,
 \baaa
&&\oo y(s)=\frac{1}{t}\int_0^sy(\t)d\t,\quad\eta=\left(\begin{array}{c}
                        \oo x(s) \\
                        y(s)\\
\oo y(s)
                      \end{array}\right)\in\R^{n+2m},\breakk\quad x(s)=\left(\begin{array}{c}
                        \oo x(s) \\
                        y(s)\\ \oo y(s)
                      \end{array}\right),
                      \eaaa
                      \baaa
&&g(x,u,s)=\left(\begin{array}{c} \oo g(x_1,...,x_n,s)\\  f(x_1,...,x_n,u,v,s) \end{array}\right),\quad\breakk
\s(x,u,s)=\left(\begin{array}{c}\oo \s(x_1,...,x_n,s)\\0_{\R^{1\times n}}\end{array}\right),
\\
\\&& h(x,u,s)=\Bigl[\Bigl(x_1-\sum_{i=1}^m u_i\Bigr)x_2\Ind_{\{0\le x_{n+i}\le C,\ i=1,...,m\}}
\breakk+\phi\left(\{x_{n+i}\}_{i=1}^m,\{x_{n+m+i}\}_{i=1}^m,t\right)
\Bigr],
 \\
 &&\Phi(x)=x_2\sum_{i=n+1}^{n+m} x_{i},
\eaaa
where  $x=(x_1,...,x_{n+2m})^\top$,
\baaa
&& f_i(x,u,v,s)= (u_{i}-\a_i x_{n+i}\breakk-\sum_{j:\ j\neq i}\b_{ij}(v_{ij}))\Ind_{\{0\le x_{n+i}\le C\}},\quad i=1,...,m,
 \\&&
  f_i(x,u,v,s)= -s^{-2} x_{n+i}+ s^{-1} x_{n+i-m},\quad\breakk i=m+1,...,2m.
 \eaaa
\section{The dynamic programming approach}\label{secHJB}
The state equations for problems  (\ref{K})  and (\ref{DKm})   are degenerate which make them difficult for analysis.
In addition, they have
 discontinuous
coefficients. To overcome this last  feature, let us approximate
the problem as the following.
Let continuously twice differentiable functions  $I_{\e,i}:\R\to [0,1]$ be defined for $\e>0$ such that
$I_{\e,i}(y)=0$ for $y\notin (0,C_i)$, $I_\e(y)=1$ for $y\in (\e,C_i-\e)$, and such that
 $I_{\e,i}(y)$ is non-decreasing in $y<\e$ and $I_{\e}(y)$ is non-increasing in $y>C_i-\e$.

 \index{Let $I_\e$ be defined similarly with $C_i$ replaced by $C$.
In the case of problem (\ref{K}), let us define functions
 \baaa
&&f_\e(x,u,t ) = (u-\a x_{n+1})I_\e(x_{n+1}),\\
&&h_\e(x,u,t ) = (x_1-u)x_2 I_\e(x_{n+1}).
\eaaa
\par
In the case of problem (\ref{DKm}),}

Let us define   $ f_\e=(f_{\e,1},...,f_{\e,2m})$
and $h_\e$ as
\baaa
&&f_{\e, i}(x,u,v,s)   = (u-\a_i x_{n+i} -\sum_{j:\ j\neq i}\b_{ij}(v_{ij}))I_\e(x_{n+i}),\\
&& h_\e(x,u,s) =
\left(x_1-\sum_{i=1}^m u_i\right)x_2\prod_{i=1}^n I_{\e,i}(x_{n+i})\breakk
+\phi\left(\{x_{n+i}\}_{i=1}^m,\{x_{n+m+i}\}_{i=1}^m,s\right).
\eaaa
Let $g_\e$ be defined similarly to $g$ with $f$ replaced by  $f_\e$.

The functions $f_\e$ and $h_\e$ approximate functions $f$ and $h$, respectively.
In addition, they are continuously differentiable.  This allows to apply the dynamic programming
approach to the  following
control  problem that approximates the original problem:
  \baa
&&\hbox{Maximize}\quad\E\left[\int_t^T h_\e(x(s),u(s),s)ds+\Phi(x(T))\right]  \quad\breakk \hbox{over}\quad
(u,v)\in \UU\nonumber \\
&&\hbox{subject to}\quad  \breakk d x(s)=g_\e(x(s),u(s),v(s),s)ds+\s(x(s),s)dw(s), \quad\breakk x(0)=x_0.
 \label{Ke}
\eaa
The corresponding optimal value function is \baa
 J_\e(x,t)
 \defi
\sup_{(u,v)\in\UU}\E\Bigl\{\int_t^T h_\e(x(s),u(s),s)ds \break+\Phi(x(T))\Bigl|
x(t)=x\Bigr\}. \label{J1}\eaa
In particular, $J_0(x,t)$ is the optimal value function  for  problem (\ref{BO1m}); by Theorem \ref{Th1M}, this is also
 the optimal value function  for problem (\ref{BKm}) and  problem (\ref{DKm}).
\par
Let $D=\R^{n+1}\times [0,T]$ for problem (\ref{K}),  and
let $D=\R^{n+2m}\times [0,T]$ for problem (\ref{DKm}).

Let $\W$ be the set of  continuous functions $w:D\to \R$  such that
there exists $c>0$ such that $|w(x,t)|\le c(|x|+1)$.
Let $\W_1=\{w\in \V:\, w'_t\in\W, \, w'_{x_i}\in \W\, (\forall i)\}$.
\begin{theorem}\label{theorem1} Assume that there are no restrictions (\ref{uvy})
(i.e.,  $\RR\equiv\R^m\times\R^{m\times m}$  therein). In this case,
For $\e>0$,
the value function $J=J_\e(x,t)$ is a solution of the following
Hamilton-Jacobi-Bellman equation
 \baa
&&J_t'+  \max_{\{u_i\},\{v_{ij}\}}  \{ J_x'g_\e +  h_\e\}+
{\scriptstyle \frac{1}{2}}{\rm Tr}\left(\b^\top J_{xx}''\b\right) = 0,\nonumber
\\
&&J(x,T)= \Phi(x). \label{BelEq}\eaa
The maximum in the HJB equation  is taken over
$u=\{u_i\}\in\R^m$, $v=\{v_{ij}\}\in\R^{m\times m}$ such that $u_i\in [-L_i,L_i] ,i=1,...,m$,
$v_{ij}\in [-\oo L_{ij},\ww L_{ij}]$, $i,j=1,...,m$, and  $(u,v)\in \RR(x_1,x_{n+1},...,x_{n+m})$.  Boundary value problem (\ref{BelEq}) has
an unique solution in the class of $\W_1$. The HJB equation holds as an equality that is
satisfied a.e. for  $(x,t)\in D$.
\end{theorem}

The following theorem establishes the way to approximate the solution $J_0$ of the original problem.
\begin{theorem}\label{ThM} The  solution $J_0$ of the original problem can be approximated as
\baa J_0(x_0,0)=\lim_{\e\to 0} J_\e(x_0,0).\label{clim2}\eaa
\end{theorem}

Hamilton-Jacobi-Bellman equation (\ref{BelEq}) can be solved via  backward
calculation after discretization and  transition to finite differences; see examples in \cite{BLM}. However, for a large $n+2m$, numerical implementation will be challenging.

 The dimension  of the HJB equation is defined by the number of factors used in the model.
For example,  if $\phi\equiv 0$, $p(t)$ is non-random, and if $S(t)$ can be modelled via a one-dimensional equation, then
 we can select $n=1$ and $m=0$. In this case, the state space will be two dimensional.
Another example:  if $m=1$ (i.e. we consider one battery only), and if $(p(t),S(t))$  can be modelled by two dimensional equation, then $n=2$. In this case, the dimension of the state space is $n+2m+1=5$. Note that modelling of energy prices and the production rate is a non-trivial task and may require rather a
 large number of factors \cite{H4}.

We address some possible ways to overcome  the problem of high dimension   in Section \ref{secD} below.

\section{Pathwise optimization for estimation of the value function}\label{secD}
Let $\oo\UU$ be the class of  all pairs $(u(t),v(t))$ of random processes with values in
$\R^m\times\R^{m\times m}$ such that (\ref{uv})-(\ref{uvy}) holds; these processes are
not necessarily $\F_t$-adapted, but their values are $\F_T$-measurable, meaning that
their path  is known for the controller at time $T$.

Let $N=m+m(m+1)/2$, and  let us associate a process  $(u,v)\in\oo\UU$ with
the $N$-dimensional  process $z_{u,v}(t)$ formed as a vector with the components  \baaa
(\{u_i(t)\}_{i=1}^m, \{v_{ij}(t)\}_{i=1,...,m, j=i+1,..,m}).\eaaa  (We excluded the components
$\{v_{ij}(t)\}_{i=1,...,m,\, j=1,..,i}$ since they are uniquely defined by $z_{u,v}(t)$ by the restrictions on $v$). 

In this section, we assume that the filtration $\{\F_t\}$ is generated by a $N$-dimensional Wiener process $w(t)$.
In particular, this case include the stochastic diffusion model described in Section \ref{SubsecDif}, given that  $n=N$.

Let $\Lambda$ be the space of all random processes $\lambda(\cdot)$ such that $\lambda(t)$ is $\F_T$-measurable  for all $t\in[0,T]$ and $\E\int_0^T\lambda(t)^2dt<+\infty$.
Let
$\Lambda$ be the  closure of the set of all processes from $\oo\Lambda$ that are $\{\F_t\}$-adapted.

For $\lambda\in \Lambda$ and $k=0,1,2,..$, let
\baaa
&&M(t)=\int_0^t \lambda(s)dW(s), \quad \mu(t)=M(T)-M(t),\quad\breakk
 \mu^{(0)}=\mu,\quad  \mu^{(k)}(t)=-\int_t^T\mu^{(k-1)}(s)ds.
 \eaaa

 Let $\oo\UU_C$ be the class of  all pairs $(u,v)\in\oo\UU$
 such that $y_i(t)\in[0,C_i]$ almost surely for all $t\in[0,T]$, $i=1,...,m$.
 Let $\UU_C$ be the class of  $\F_t$-adapted pairs $(u,v)\in \oo\UU_C$.

 For $(u,v)\in \oo\UU_C$, let us define
 \baa
&&F(u,v)=\E\int_t^T\Bigl[\left(p(s)-\sum_{i=1}^m u_i(s)\right)S(s)
\breakk +\phi(y(s),\oo y(s),s)\Bigr] ds +\E S(T)\sum_{i=1}^m y_i(T),
\label{Fuv}\eaa
where $y$ and $\oo y$ are such as  defined in (\ref{BO1m}).

Let $k\in\{0,1,2,...\}$ be selected, and let
\baaa \L(u,v,\lambda)=\E F(u,v)+\E\int_0^T \mu^{(k)}(t)^\top z_{u,v}(t)dt. \eaaa

 \begin{theorem}\label{ThD} Assume that
 the function $\phi(y,\oo y, t)$ is concave in $(y,\oo y)$, Then
\baa
 \sup_{(u,v)\in {\mathscr{U}}_C} \E F(u,v) =
  \inf_{\lambda\in \Lambda} \sup_{( u, v)\in \oo{\mathscr{U}}_C} \L(u,v,\lambda).
\label{infsup} \eaa
\end{theorem}

It can be noted that the assumption on  $b_{ij}=0$ above means  that the link losses for transmission of energy between batteries are not taken into account.
\par
In Theorem \ref{ThD}, $\L$ is an analog of the Lagrangian for the control problem with special  constraints that
$(u,v)\in {\mathscr{U}}_C$.
\par
 Theorem \ref{ThD}  allows to substitute stochastic  optimization over the set of adapted  controls by the set
 of non-adapted controls (i.e. controls using the information about the future). This allows to estimate the value function by
 Monte-Carlo method via simulating $(u,v,\lambda)\in \oo\UU\times\Lambda$ using approach
\cite{BS,Br,Davis,HK,Ro,D17}.  Therefore, Theorem \ref{ThD}
gives an opportunity to estimate how far from optimal is the performance of a particular  strategy.

It has to be clarified that an existence
of a saddle point does not follow from
Theorem \ref{ThD};  also,  this theorem does not give a way to derive an optimal strategy $u$.

\section{Proofs}\label{secP}
\par
 {\it Proof of Theorem \ref{ThD1}.} Let
 $(y(\cdot),u(\cdot))$ be an admissible process  for problem (\ref{O1}). It is  also an admissible process
  for problem (\ref{K1}), and the values of the quality criterions for $(y(\cdot),u(\cdot))$  are the same for both problems. Hence the optimal value (i.e., the value of the expectation in the performance criterion  for problem (\ref{O1}) is less or equal the the optimal value  for problem (\ref{K1}). Further, let
  $(y(\cdot),\ww u(\cdot))$ be an admissible process  for problem (\ref{K1}). The equation for $y$  in (\ref{K1}) is such that
  $y(s)\in[0,C]$ for all $t$. In addition, the same $y$ will be generated by the control  $\w u(t)$ such that
  $\w u(t)=u(t)$ if $y(t)\in (0,L)$ and $\w u(t)= \a y(t)$ if $y(t)=0$ or $y(t)=L$. It follows that the process
  $(y(t),\w u(t))$ is also admissible for problem (\ref{O1}). Furthermore,  the control $u$  will not
  perform better than $\w u$ because of the presence of the indicator function under the integral in (\ref{K1}).  Hence the optimal value (i.e., the value of the expectation in the performance criterion) for problem (\ref{K1}) is less or equal the the optimal value  for problem (\ref{O1}).
This proves  Theorem \ref{ThD1}. $\Box$

The proof of Theorem \ref{Th1M} is similar.

\par
 {\it Proof of Theorem \ref{theorem1}.}
 Theorem 4.1.1 \cite{Krylov}, p.165,  implies that  $J=J_\e$ satisfies the
 HJB equation such that it has unique solution in
$\W_1$ and all components of $J'_x$  belong to $\W$. This HJB equation
holds  in the sense of an equality of the distributions.  Theorem 4.4.3 \cite{Krylov}, p.192, implies that  $J'_t\in \V$ as well.  $\Box$.
 \par
 {\it Proof of Theorem \ref{ThM}.}
 Let $(x,t)$ be fixed, and let
\baaa
F_\e(u,v)=\E\Bigl\{\int_t^T h_\e(x(s),u(s),s)ds +\Phi(x(T))\Bigr),
\eaaa
 where $(u,v)\in\UU$, and where  $x$ is the corresponding solution of the differential equation in (\ref{DKm}).

 Let $\d>0$, and let $(\w u,\w v)$ be such that $F_0(\w u,\w v))\ge J_0(x,t)-\d/2$. Let   $\w y(t)$ be the corresponding vector process. Let $\w u_\e=\{\w u_{\e,i}\}$ and  $\w v_\e=\{\w v_{\e,ij}\}$  be such that
 $\w u_{\e,i}(t)=\w u_i(t)$ and $\w v_{\e,ij}(t)=\w v_{ij}(t)$ if $\w y_i(t)\in (\e,C-\e)$, and where
 $\w u_{\e,i}(t)=\a_i \w y_i(t)$ and $\w v_{\e,ij}(t)=\w v_{\e,ji}(t)=0$     if $\w y_i(t)\in \{\e,C-\e\}$. We have that there exists $\e>0$ such that
  $F(\w u_\e,\w v_\e)\ge F(\w u,\w v)-\d/2$. Therefore, for an arbitrarily small  $\d>0$, one can find
  $(\w u_\e,\w v_{\e})$ such that
  $F(\w u_\e,\w v_\e)\ge F(\w u,\w v)-\d\ge J_0(x,t)-\d/2$.  Then the proof follows. $\Box$
 \par
 {\em Proof of Theorem \ref{ThD}}.
 It can be shown that \index{$M(T)\in\LL_2(\O,\F_T,\P;\R^N)$ and}  $M(t)=\E\{M(T)|\F_t\}$ (i.e., this is the conditional expectation given $\F_t$).
 Similarly to the  proof of Theorem 4.1  \cite{D17},
 we obtain   that  $(u,v)\in \oo{\mathscr{U}}_C$
 belongs to $\mathscr{U}_C$
if and only if
 \baaa
 \E\int_0^T \mu(t)^\top z^{(k)}(t) dt=\E\int_0^T \mu^{(k)}(t)^\top z(t)dt=0
 \eaaa  for any $k\ge 0$ and any $\lambda\in \Lambda$.
 Here   $z^{(0)}=z_{u,v}$,
 $z^{(k)}(t)=\int_0^tz^{(k-1)}(s)ds$, $k=1,2,3,...$. This is why
 $\L(u,v,\lambda)$ can be used as a Lagrangian for the problem with the constraint
 that $(u,v)\in{\mathscr{U}}_C$.

Further, the set  $\oo{\mathscr{U}}_C$ is
convex and closed in the square integral metric. By the assumtions on $b_{ij}$,
we have that $y$ and $\oo y$ depend linearly on $(u,v)\in\oo{\mathscr{U}}_C$. By the assumptions on $\phi$,   the function $F(u,v)$ is
concave on $\oo{\mathscr{U}}_C$. Then the proof of Theorem \ref{ThD}
follows  The statement of the theorem follows Propositions 1.2 and 2.3 \cite{Ekland}, Chapter 4,
similarly  to the  proof of Theorem 4.1  \cite{D17}.
 $\Box$

\section{Conclusion}\label{secC}

The paper suggests a compact and yet comprehensive model for decision making under uncertainty
for a  of small or medium
producer of energy operated a plant and a system of batteries.
The problem setting  takes into account the topology of the system of batteries and losses for transferring and storage of energy.
A method of calculating an optimal operating algorithm for storing and dispatching strategy is suggested. The optimality is defined by a performance criterion that
takes into account  the expected monetary return for the producer.
In addition, the performance criterion takes into account preferable regimes for batteries
charging and discharging that may help to  prolong  the battery life.

The model includes two processes that cannot be controlled by the operator: the  production rate  $p(t)$ and the energy price  $S(t)$. If these processes are predictable, then calculation  of the optimal strategy requires to solve a  multi-dimensional deterministic optimal control problem. The paper is focused on the setting
where processes $p(t)$ and $S(t)$ are random, currently observable, and unpredictable, with known stochastic evolution law.   This leads to a multi-dimensional stochastic optimal control problem.
A method of  solution is developed for  the model where $(p(t),S(t))$
is  part of stochastic Markov diffusion process.
  There are no restrictions on the  choice of the parameters of this diffusion process;  in this sense, the model is very general.
In particular, this methods is also applicable for the case where  the process  $(p(t),S(t))$ is non-random.

The Markov  diffusion  models are quite common in engineering and economics; their
main restriction is that they use continuous processes without  jumps.
The model introduced in this paper is quite flexible and allow many modifications that were not
executed in the present paper to avoid overloading by technical details.

 Instead of the diffusion Markov model described in Section \ref{secM},
many other stochastic models  can be accommodated for the dynamics of  $(p(t),S(t))$.
For example, the dynamic programming  method can be extended on the jump diffusion model using the approach \cite{OS}.  A case of general right-continuous processes $(p,S)$ could  be considered similarly to \cite{BD2} The continuous time model could be  replaced by a  discrete time model.
We leave this for the future research.

\index{Also, a case of general right-continuous processes $(p,S)$ can be considered similarly to \cite{BD1,BD2}; this
would require more complicated analysis. It can be also noted that the continuous time model could be  replaced by a Markov discrete time model.}

Potentially,  uncertainty can be taken into account with stochastic models
replaced by interval type uncertainty.  We leave this for the future research as well.

In the present form, the   model is based on the  continuous time processes. However,
this is rather technical assumption since a similar discrete
time model can be obtained via straightforward discretization.

Numerical implementation of the suggested method  would require time discretization
 and  backward solution of the dynamic programming
equation. However, the state space dimension for this equation
will  be high for a large number of factors arising in a multi-battery setting.
This means that numerical solution could  require significant computational efforts.
To partially  address this, we suggest to estimate  of the optimal value function
using the pathwise optimization approach.
We leave further development of the numerical methods for the future research.

Usually, decision making in the stochastic framework relies on forecasting  the underlying processes. In  stochastic evolution models  such as  (\ref{diff})-(\ref{oox}), the forecasting is assumed in an implicit form:
 a model implies  forecasts for the underlying  processes, for example, the conditional expectations of the future values given the observations. These
  forecasts and the corresponding  errors
 can be calculated analytically for a particular model such as  (\ref{diff})-(\ref{oox}). One would need to apply statistical inference and  forecasting methods  to
  calibrate the parameters for  a model. However, this requires special consideration beyond the scope of this paper. We leave  this  for the future research.

Several other important factors have been left for the future research. 
Possible cooperation with the grid operators and other
producers could be included in the model.  Some restrictions on batteries could be imposed, for example,
on simultaneous charging or discharging.
It  would be interesting to consider
the problem of optimal selection of the topology of the system of batteries to minimize the losses. 


 \index{\subsubsection*{Acknowledgment}
  This work was
supported \comm{financially supported}  by Grant 08-08 of Government of Russian Federation. \blind{ARC} grant of \blind{Australia DP120100928}  to the author.}

\end{document}